\documentclass[12pt,a4paper]{article}

\usepackage{graphicx}
\usepackage{subfig}
\usepackage{amsmath,amsthm,amssymb}
\usepackage{esint}
\usepackage{mathrsfs}
\usepackage{appendix}
\usepackage{hyperref}
\newcommand{\RomanNumeralCaps}[1]

\usepackage[left=2.2cm,right=2.2cm,top=3cm,bottom=3.5cm]{geometry}

\usepackage{color}
\usepackage{array}
\usepackage{url}
\usepackage{changepage}
\usepackage{float}

\DeclareMathOperator{\bnabla}{\boldsymbol{\nabla}}

\newtheorem*{theorem*}{Theorem}

\theoremstyle{definition}

\theoremstyle{remark}

\usepackage{authblk}

\usepackage{mathtools}
\usepackage{amssymb}
\usepackage{appendix}

\numberwithin{equation}{section}

\usepackage[sort&compress]{natbib}

\newcommand{\bom}{\mbox{\boldmath$\omega$}}

\newcommand{\bu}{\mbox{\boldmath$u$}}

\newcommand{\bx}{\mbox{\boldmath$x$}}

\newcommand{\bB}{\mbox{\boldmath$B$}}

\newcommand{\I}{\int_{V}}
\newcommand{\Div}{\mbox{div}\,}
\newcommand{\divu}{\mbox{div}\,\bu}
\newcommand{\capH}{H}
\newcommand{\capHC}{H_{c}}
\newcommand{\bdJ}{\mbox{\boldmath$J$}}
\newcommand{\bel}{\begin{equation}\label}
\newcommand{\ee}{\end{equation}}
\newcommand{\beq}{\begin{eqnarray}\label} 
\newcommand{\eeq}{\end{eqnarray}} 
\newcommand{\bc}{\begin{center}} 
\newcommand{\ec}{\end{center}} 
\newcommand{\ben}{\begin{enumerate}}
\newcommand{\een}{\end{enumerate}}
\newcommand{\bit}{\begin{itemize}}
\newcommand{\eit}{\end{itemize}}
\newcommand\shalf{\ensuremath{{\scriptstyle\frac{1}{2}}}}

\usepackage[T1]{fontenc}
\usepackage[polish,english]{babel}
\usepackage[utf8]{inputenc}

\setcitestyle{aysep={}}

\allowdisplaybreaks

\begin{document}

\title{\large\bf On the removal of the barotropic condition in helicity studies of the compressible Euler and ideal compressible MHD equations}

\author[]{Daniel W. Boutros\footnote{Department of Applied Mathematics and Theoretical Physics, University of Cambridge, Cambridge CB3 0WA UK. Email: \textsf{dwb42@cam.ac.uk}} \space and John D. Gibbon\footnote{Department of Mathematics, Imperial College London, London SW7 2AZ, UK. Email: \textsf{j.d.gibbon@ic.ac.uk}}}

\date{January 29, 2025}

\maketitle

\begin{abstract}
The helicity is a topological conserved quantity of the Euler equations which imposes significant constraints on the dynamics of vortex lines. In the compressible setting the conservation law only holds under the assumption that the pressure is barotropic. We show that by introducing a new definition of helicity density $h_{\rho}=(\rho\bu)\cdot\mbox{curl}\,(\rho\bu)$ this assumption on the pressure can be removed, although $\I h_{\rho}dV$ is no longer conserved. However, we show for the non-barotropic compressible Euler equations that the new helicity density $h_{\rho}$ obeys an entropy-type relation (in the sense of hyperbolic conservation laws) whose flux $\bdJ_{\rho}$ contains all the pressure terms and whose source involves the potential vorticity $q = \bom\cdot \nabla \rho$. Therefore the rate of change of $\I h_{\rho}dV$ no longer depends on the pressure and is easier to analyse, as it only depends on the potential vorticity and kinetic energy as well as $\divu$. This result also carries over to the inhomogeneous incompressible Euler equations for which the potential vorticity $q$ is a material constant. Therefore $q$ is bounded by its initial value $q_{0}=q(\bx,\,0)$, which enables us to define an inverse resolution length scale $\lambda_{H}^{-1}$ whose upper bound is found to be proportional to $\|q_{0}\|_{\infty}^{2/7}$. In a similar manner, we also introduce a new cross-helicity density for the ideal non-barotropic magnetohydrodynamic (MHD) equations.
\end{abstract}


\vspace{0.1cm} \noindent \textbf{Mathematics Subject Classification:} 76N99 (primary), 76W05, 76B99 (secondary)

\section{\large Introduction}\label{sect1}

It was first shown by \citet{Helm1858} that for an ideal, barotropic fluid with conservative body forces, vortex lines are transported by the flow. \citet{Kelvin1868} (later Lord Kelvin) then recognised that any knots and linkages in these lines are conserved. Almost a century later, \citet{Moreau1961} showed that the helicity $\mathcal{H}$ is conserved, a result that was later proved independently by \citet{Moffatt1969}, who went further in recognizing its relationship with the magnetic helicity invariant of \citet{Woltjer1958}, and with the cross-helicity invariant $\mathcal{H}_{c}$ of the ideal magnetohydrodynamic (MHD) equations \citep{Woltjer1958b}. 
\par\smallskip%\noindent
In a fluid flow with velocity vector $\bu$ and vorticity vector $\bom = \nabla \times \bu$ the helicity is defined to be a volume integral of the form $\mathcal{H}=\I\bu\cdot\bom\,dV$.  In this context, \citet{Moffatt1969} showed explicitly that it is also a measure of the degree of vortex line linkages in a localised disturbance. This reference \citep{Moffatt1969} has been the foundation for a significant body of work on helicity applied to knots and linkages in both ideal fluids and MHD\,: see also \citet{Pouquet1976}, \citet{Berger1984}, \citet{Moffatt1978,Moffatt1985,Moffatt1990} and \citet{MoffattRicca1992}. Moreover, it was shown in \citet{Enciso2016} that the helicity is the only integral invariant for a general volume-preserving flow. We refer to the paper by \citet{Moffatt2014} for an overview of results and a wide range of references in the subject. The effect of helicity dynamics on the energy cascade, as well as the energy dissipation rate, has been investigated in \citet{Biferale2012,Capocci2023,Linkmann2018} (see references therein), see also \citet{Biferale2013}. In fact, tracking the helicity density has now become a standard diagnostic tool in the study of vortical structures in large scale numerical simulations of both incompressible Euler and Navier-Stokes flows\,: see \citet{Kerr2018,Kerr2023}. For more general recent references in incompressible Navier-Stokes turbulence, see \citet{ISY2019}, \citet{PK2020} and \citet{BPB2022}.
\par\smallskip%\noindent
In the setting of ideal compressible flows, a (nearly) universal feature in the study of helicity dynamics has been the assumption of barotropicity of the fluid\,: that is, the pressure $P$ is taken to be solely a function of the fluid density $\rho$, so iso-surfaces of pressure and density are parallel. As has been described in \citet{Thorpe2003}, the results of Helmholtz and Kelvin on vortex lines and structures were introduced to the meteorological and oceanographic communities by \citet{VB1898} (in which rotational effects were included) through his circulation theorem\footnote{This is also called the Poincar\'e-Bjerknes circulation theorem \citep{Poincare1893}.}. The essential role played by barotropicity had been pointed out earlier by \citet{Silber1896}. In meteorological reality, the validity of the barotropic assumption is somewhat limited. For instance, regions of barotropicity in the atmosphere have a uniform temperature distribution and are distinguished by the absence of fronts, so the barotropic assumption is generally restricted to the tropics \citep{RN2013,Vallis2019,IR2021}, whereas in the mid-latitudes the atmosphere is generally baroclinic. For instance, the Unified Model run by the UK Meteorological Office assumes that the pressure is a function of both density and temperature \citep{UM2005}. The barotropic assumption has greater validity, however, in the study of stellar interiors or of the interstellar medium \citep{Yokoi2013}. One common class of barotropic models used in astrophysics are polytropic fluids, where it is assumed that $P \propto \rho^{\gamma}$ for $\gamma \geq 1$.
\par\smallskip%\noindent
In the non-barotropic case, a generalised helicity (which is a conserved quantity) has been introduced in \citet{Mobbs1981}\,: see also \citet{Gaffet1985} and \citet{Salmon1988}.  Similarly, a generalised cross helicity can be introduced for the non-barotropic MHD equations \citep{Yahalom2017}. As noted in \citet{Webb2018} (see page 183), these generalised helicities are nonlocal quantities as they depend on a nonlocal variable, namely the Lagrangian time integral of the temperature. The goal of this paper is to introduce a generalised helicity which solely depends on {\em local} variables, even though it might not be fully conserved. Due to its local nature, the dynamics of this new generalised helicity is easier to interpret.
\par\smallskip%\noindent
In particular, the aim of this paper is to explore the circumstances under which the barotropic condition can be removed and to study the consequences of its removal. In order to explain briefly how \citet{Moreau1961}  and  \citet{Moffatt1969} used the barotropic assumption, let us consider an ideal compressible fluid with velocity vector $\bu$, density $\rho$, pressure $P$ and vorticity $\bom = \mbox{curl}\,\bu$. In their standard form, without the inclusion of the (specific) internal energy (see \S\ref{sect2b}), the barotropic compressible Euler equations take the form
\bel{comp1a}
\rho\frac{\mathrm{D}\bu}{\mathrm{D}t} =   - \bnabla P\,, \qquad\qquad 
\frac{\mathrm{D}\rho}{\mathrm{D}t} + \rho\divu = 0\,,
\ee
where $\mathrm{D}/\mathrm{D}t = \partial_{t}+\bu\cdot\bnabla$ is the material derivative and $P$ is solely a function of $\rho$. The equation for the vorticity is
\bel{comp1b}
\frac{\mathrm{D}\bom}{\mathrm{D}t} + \bom\divu = \bom\cdot\bnabla\bu 
- \bnabla\left(\rho^{-1}\right)\times\bnabla P\,.
\ee
We note that the cross product vanishes as $\bnabla P$ and $\bnabla\left(\rho^{-1}\right)$ are both parallel to $\bnabla \rho$, due to the barotropic assumption. Therefore the evolution of the helicity density can be written as an entropy-type relation (in the sense of hyperbolic conservation laws) of the form
\bel{comp2}
\partial_{t}h + \mbox{div}\, \bdJ_{\pi} = 0\,,
\ee
where $h= \bu\cdot\bom$, $\Pi (\rho) = \int_{0}^{\rho}\eta^{-1}P' (\eta) d\eta$ and
\bel{comp3}
\bdJ_{\pi} = h\bu + \bom\left(\Pi - \shalf|\bu|^{2}\right)\,.
\ee
Thus, with suitable boundary conditions on the domain such as periodic boundary conditions, it follows that $\mathcal{H} = \I \bu\cdot\bom\,dV$ is conserved. In fact, \citet{Moffatt1969} showed that the time derivative of the helicity integrated over a moving domain $\mathcal{V}(t)$ transported by the fluid is equal to a perfect divergence. Therefore its evolution on such a moving domain is only determined by boundary terms.
\par\smallskip%\noindent
The above calculation illustrates the point that if the barotropic assumption is to be dropped then the helicity $\mathcal{H}=\I h\,dV$ is no longer conserved, because in equation \eqref{comp2} the term $\mbox{div}  \left( \Pi \bom \right)$ is replaced by $\rho^{-1}\mbox{div} (\bom P)$ and the additional term $\bu \cdot (\bnabla (\rho^{-1}) \times \bnabla P)$ appears. Therefore a change in the definition of the helicity $\mathcal{H}$ is required.  It is hardly surprising that without the barotropic assumption, the property that vortex lines are transported with the flow will be lost. We introduce a new definition of helicity density $h_{\rho}$ and helicity $H$ as
\bel{Hdef}
h_{\rho} = (\rho\bu)\cdot\mbox{curl}\,(\rho\bu)=(\rho\bu)\cdot(\rho\bom)\qquad\mbox{with}\qquad \capH = \I h_{\rho}\,dV\,.
\ee
In the context of helicity dynamics without the barotropic assumption, the challenge is to estimate the growth or decay of a topological quantity without explicitly imposing any assumptions on $P$. In \S\ref{sect2} of this paper we will show that $h_{\rho}$ obeys an entropy-type law
\bel{cl1}
\partial_{t}h_{\rho} + \mbox{div}\bdJ_{\rho} = \sigma_{\rho}\,.
\ee
The pressure $P$ appears only in $\bdJ_{\rho}$ and thus disappears under integration over a periodic domain. In the case of the inhomogeneous incompressible Euler equations, a potential physical interpretation of these formal results lies in using the helicity $H$ and the constant total energy $E_{0}$ to define an inverse length scale. In analogy with the Kolmogorov length scale, in \S\ref{sect2a} we define this inverse length scale based on the average rate of change of the helicity as follows
\bel{lamdefA}
\left[\lambda_{H}\right]^{-1} = \left(\frac{\left\langle \left\lvert \frac{\mathrm{d}\capH}{\mathrm{d}t} \right\rvert \right\rangle}{\varrho_{0}^{1/2} E_{0}^{3/2}}\right)^{2/7}\,, \qquad\qquad 
\varrho_{0} = L^{-3}\I \rho dV\,.
\ee
The brackets $\langle \cdot \rangle \coloneqq \frac{1}{T} \int_0^T \cdot \, dt$ denote a finite-time average, $\varrho_{0}$ denotes the domain averaged density and the energy is defined by
\bel{Edef}
E_{0} = \I \mathcal{E}_{0}dV\qquad \mbox{and}\qquad \mathcal{E}_{0}=\shalf\rho|\bu|^{2}\,.
\ee
The length scale $\lambda_{H}$ could be interpreted as the smallest length scale on which there are significant variations of $H$ and hence significant topological variations. The source term $\sigma_{\rho}$ is twice the product of the potential vorticity $q =\bom\cdot\bnabla\rho$ and the energy density. It is shown in \S\ref{sect2a} that this leads to the estimate
\bel{lamest1A}
\left[\lambda_{H}\right]^{-1} \leq \left( \frac{4}{E_0 \rho_0} \right)^{1/7} \lVert q_0 \rVert_\infty^{2/7} \,,
\ee
where $\|q_0\|_{\infty}$ is the spatial maximum norm of the (initial) potential vorticity, which is a material constant under the dynamics. We emphasise that for the inhomogeneous incompressible Euler equations the pressure solves an elliptic problem and the barotropic assumption is not relevant in this context, as the pressure also depends on velocity derivatives. As the canonical helicity is not a conserved quantity for these equations, the modified helicity $H$ provides a new topological quantity whose dynamics is easy to interpret and its growth is bounded. 
\par\smallskip%\noindent
A similar calculation for the fully compressible Euler equations (including the specific internal energy) is explained in \S\ref{sect2b} with the complication that $\sigma_{\rho}$ has an extra term involving $\divu$ that weights the helicity density $h_{\rho}$. Finally in \S\ref{MHD} we introduce a non-barotropic cross-helicity density $h_{c}=\rho\bu\cdot\bB$ for the ideal compressible MHD equations.

%%%%%%%%%%%%%%%%%%%%%%%%%
\section{\large The evolution of $h_{\rho}$ for both the inhomogeneous incompressible and the compressible $3D$ Euler equations}\label{sect2}

It is known that strong solutions of the $3D$ incompressible Euler equations can develop singularities in a finite time \citep{Elgindi2021}\,: see also \citet{Drivas2023}  for a recent survey and references therein. In the compressible case it is also known that smooth solutions blow up in finite time \citep{Sideris1985}. Once singularities or shocks develop, solutions could be too irregular to perform the manipulations needed to obtain the results in this paper\,: for example, the vorticity $\bom$ might no longer be a well-defined pointwise quantity. Thus it should be understood that the results in the following sections are only valid for time intervals when sufficiently smooth solutions exist for both the Euler and ideal MHD equations. Sufficient regularity conditions for solutions of the incompressible Euler equations to conserve the helicity, can be found in \citet{CCFS2008} and \citet{BT2024} and references therein.

%%%%%%%%%%%%%
\subsection{\small Results for the inhomogeneous incompressible $3D$ Euler equations}\label{sect2a}

We recall that the $3D$ inhomogeneous incompressible Euler equations are given by
\bel{iie1}
\rho\frac{\mathrm{D}\bu}{\mathrm{D}t} =   - \bnabla P\,, \qquad \divu = 0\,,\qquad 
\frac{\mathrm{D}\rho}{\mathrm{D}t} = 0\,.
\ee
In this model the density $\rho$ is allowed to vary but $\bu$ remains divergence-free. The divergence-free condition implies that $P$ must satisfy an elliptic equation that involves derivatives of the velocity
\begin{equation}
\bnabla \cdot \left( \frac{1}{\rho} \bnabla P \right) = - (\bnabla\otimes\bnabla) : (\bu\otimes\bu) \,,
\end{equation}
and so an imposition of barotropicity is invalid. The vorticity satisfies the following equation
\bel{inhomegvorticity}
\frac{\mathrm{D}\bom}{\mathrm{D}t} = \bom\cdot\bnabla\bu 
- \bnabla\left(\rho^{-1}\right)\times\bnabla P\,.
\ee 
We note that the canonical helicity $\mathcal{H}$ is not conserved by solutions of the inhomogeneous incompressible Euler equations, as the integral $\I \bu \cdot \big[\bnabla\left(\rho^{-1}\right)\times\bnabla P \big] \,dV$ generally does not vanish. The new helicity $H$ introduced in (\ref{Hdef}) is also not a conserved quantity, but we will find that its growth is bounded and its evolution depends only on local quantities. First we need to recall two basic identities.
\par\smallskip%\noindent
Firstly, it is not difficult to show that $\rho$, $\bu$ and $P$ obey the conservation law
\bel{fr1}
\partial_{t}\mathcal{E}_{0} + \mbox{div}\left\{\left(\mathcal{E}_{0}+P\right)\bu\right\}= 0\,,
\ee
where $\mathcal{E}_{0} = \shalf\rho|\bu|^{2}$ is the energy density. Therefore on a periodic domain the energy $E_{0} = \int_{V}\mathcal{E}_{0}\,dV$ is conserved. 
\par\smallskip%\noindent
Secondly, we recall Ertel's theorem \citep{Ertel1942} for the potential vorticity $q = \bom\cdot\bnabla\rho$. Using the evolution equation for $\bom$ in (\ref{inhomegvorticity}), it can be shown that $q$ satisfies
\bel{fr2}
\frac{\mathrm{D}q}{\mathrm{D}t} = \bom\cdot\bnabla\left(\frac{\mathrm{D}\rho}{\mathrm{D}t}\right) - \left[\bnabla(\rho^{-1})\times\bnabla P\right]\cdot\bnabla\rho = 0\,,
\ee
where terms of the type $(\bnabla\rho)\cdot(\bom\cdot\bnabla\bu)$ cancel, in tandem with the vanishing of the pressure term.  It therefore follows that $q$ is a material constant. 
\par\smallskip%\noindent
The material derivatives of $\rho\bu$ and $\rho\bom$ are easily found, thereby giving 
\beq{h1}
\frac{\mathrm{D}h_{\rho}}{\mathrm{D}t} &=& \rho^{2}\bu\cdot(\bom\cdot\bnabla\bu) 
-  \left\{\rho\,\bom\cdot\bnabla P - \bu \cdot (\bnabla\rho\times\bnabla P)\right\}\,.
\eeq
We note that the terms involving the pressure in (\ref{h1}) form a perfect divergence
\beq{h1A}
\rho\,\bom\cdot\bnabla P - \bu\cdot[\bnabla\rho\times\bnabla P] &=& \bnabla P\cdot\left[ \rho\,\bom +\bnabla\rho\times\bu\right]\nonumber\\
&=& \bnabla P \cdot \mbox{curl}\,(\rho\bu)\nonumber\\
&=& \mbox{div}\left\{P\mbox{curl}\,(\rho\bu)\right\}\,.
\eeq
Equation (\ref{h1}) then becomes
\beq{eb3}
\partial_{t} h_{\rho} + \mbox{div}\,\{h_{\rho}\bu + P\mbox{curl}\,(\rho\bu)\} &=& \rho^{2}\bom\cdot\bnabla(\shalf|\bu|^{2})\nonumber\\
&=& \bom\cdot\bnabla\left(\shalf\rho^{2}|\bu|^{2}\right) - |\bu|^{2}\bom\cdot\bnabla(\shalf\rho^{2})\nonumber\\
&=& \mbox{div}\left\{\shalf\bom\rho^{2}|\bu|^{2}\right\} - q\left(\rho|\bu|^{2}\right)\,,
\eeq
which can be written as an entropy-type relation
\bel{lem1}
\partial_{t} h_{\rho} + {\rm div}\bdJ_{\rho} = \sigma_{\rho}\,,
\ee
with the flux vector $\bdJ_{\rho}$ and the scalar source term $\sigma_{\rho}$ defined as 
\bel{Jhdef}
\bdJ_{\rho} = h_{\rho}\bu + P\mbox{curl}\,(\rho\bu)-\shalf\bom\rho^{2}|\bu|^{2}\,,\qquad \qquad
\sigma_{\rho} = - q\rho|\bu|^{2}\,.
\ee
One can therefore infer that the sign of the potential vorticity impacts whether $h_{\rho}$ increases or decreases. After integration of equation \eqref{lem1} over a periodic domain, the term ${\rm div}\bdJ_{\rho}$ disappears and one finds
\bel{lem1B}
\frac{\mathrm{d}\capH}{\mathrm{d}t} = - 2\int_{V}q\mathcal{E}_{0}\,.
\ee
The property that $q$ is a material constant means that the growth or decay of $H$ is bounded. In fact, equation \eqref{fr2} implies that $\lVert q (\cdot, t) \rVert_\infty \leq \lVert q_0 \rVert_\infty$. This immediately implies the bound
\bel{growthbound}
\left\lvert \frac{\mathrm{d}\capH}{\mathrm{d}t} \right\rvert \leq 2 \lVert q_0 \rVert_\infty E_0\,,
\ee
where $\|q\|_{\infty}$ is the maximum spatial norm of the potential vorticity. In turn, this bound constrains the globally averaged alignment between $\bu$ and $\bom$ and hence the topological dynamics. 
\par\smallskip%\noindent
Next we notice that the evolution of $H$ induces a length scale $\lambda_{H}$ of the following form
\bel{lamdefB}
\left[\lambda_{H}\right]^{-1} = \left(\frac{\left\langle \left\lvert \frac{\mathrm{d}\capH}{\mathrm{d}t} \right\rvert \right\rangle}{\varrho_{0}^{1/2} E_{0}^{3/2}}\right)^{2/7}\,,
\ee
where $\varrho_{0}$ is the average density $\varrho_{0}= L^{-3}\I \rho\,dV$. The quantity $\lambda_{H}$ can be interpreted as a resolution length scale\,: for example the smallest length scale on which there is significant topological dynamics. We can now use equation \eqref{growthbound} to find the following upper bound on this inverse length scale
\bel{lamest1B}
\lambda_{H}^{-1} \leq \left( \frac{4}{E_0 \rho_0} \right)^{1/7} \lVert q_0 \rVert_\infty^{2/7} \,.
\ee
As has been said before, $\lambda_H$ could be viewed as a cutoff length scale for the helicity dynamics.

%%%%%%%%%%%%%%%%%%%%%%%%%

\subsection{\small Results for the fully compressible $3D$ Euler equations}\label{sect2b}

The fully compressible Euler equations, including the specific internal energy $e$ (per unit mass), require an equation for $e$ in addition to those given in (\ref{comp1a}). They are given by
\bel{ce1a}
\rho\frac{\mathrm{D}\bu}{\mathrm{D}t} =   - \bnabla P\,,\qquad\mbox\qquad 
\frac{\mathrm{D}\rho}{\mathrm{D}t} + \rho\,\divu = 0\,,\qquad\qquad \rho\frac{\mathrm{D}e}{\mathrm{D}t} = - P\divu\,.
\ee
Before considering the helicity density $h_{\rho}$, let us consider the well-known formula for the full energy density 
\bel{ed1}
\mathcal{E} = \rho\left(\shalf|\bu|^{2} + e\right)\,.
\ee
It is not difficult to show that $\mathcal{E}$ satisfies the exact continuity equation (\ref{fr1}) (with $\mathcal{E}_{0}$ replaced by $\mathcal{E}$). Thus we find that the total energy $E = \I \mathcal{E}dV$ is constant for any equation of state. For the system \eqref{ce1a} to be fully determined, an equation of state for $P$ in terms of $e$ and $\rho$ is required. Our results however are independent of the choice of equation of state and therefore we do not fix a choice. We note that similar results as described below can also be obtained for different formulations of the compressible Euler equations involving temperature or entropy dynamics. This is because our results only rely on the form of the density and velocity equations. In addition, we note that the potential vorticity $q$ is not a material constant for these equations. However, it was shown in \citet{Gibbon2012} that $q$ satisfies 
\bel{GHq}
\partial_t q + \Div (q \bu) + \Div \left[ \bom \rho \divu \right] = 0\,.
\ee
One can deduce from equation (\ref{ce1a}) that $\rho\bu$ and $\rho\bom$ evolve according to
\beq{ce2}
\frac{\mathrm{D}(\rho\bu)}{\mathrm{D}t} + (\rho\bu)\divu &=& -\bnabla P \,,\\
\frac{\mathrm{D}(\rho\bom)}{\mathrm{D}t} + 2(\rho\bom)\divu &=& \rho\bom\cdot\bnabla\bu + \rho^{-1}\bnabla\rho\times \bnabla P \,.\label{ce3}
\eeq
Following the grouping of the pressure terms as in (\ref{h1A}), from (\ref{ce2}) and (\ref{ce3}) we deduce that
\bel{ce5}
\partial_{t} h_{\rho} + 2h_{\rho} \divu + \mbox{div}\left\{h_{\rho}\bu + P\mbox{curl}\,(\rho\bu)-\shalf\bom\rho^{2}|\bu|^{2}\right\} = - q\rho|\bu|^{2}\,.
\ee
With $\bdJ_{\rho}$ and $\sigma_{\rho}$ defined in (\ref{Jhdef}) and $\tilde{\sigma}_{\rho} = \sigma_{\rho} - 2h_{\rho} \,\divu$, then we find a similar entropy-type relation to (\ref{lem1})  
\bel{lem2}
\partial_{t} h_{\rho} + {\rm div}\bdJ_{\rho} = \tilde{\sigma}_{\rho} \,.
\ee
(\ref{lem2}) integrates to 
\bel{lem2B}
\frac{\mathrm{d}\capH}{\mathrm{d}t} + 2 \int_{V}h_{\rho}\,\divu\,dV = - 2 \I q \mathcal{E}_{0} dV\,.
\ee
From this equation one can observe that the helicity $H$ increases in regions of compression, while it decreases in regions of dilatation. Moreover, one can deduce an inequality of the following form
\bel{compressiblegrowthbound}
\left\lvert \frac{\mathrm{d}\capH}{\mathrm{d}t} \right\rvert \leq 2 \lVert q (\cdot, t) \rVert_\infty E_0 + 2 \lVert h_{\rho} (\cdot, t) \rVert_\infty \I |\divu (\cdot, t)|\,dV\,,
\ee
where we observe that bounds of this type are particularly useful in the perturbative regime of slightly compressible flows where the last term on the right hand side could be small.

%%%%%%%%%%%%%%%%%%%%%%%%%%%%%%%%%%%%%%%%%%%%%%%%%%
\section{\large Cross-helicity in ideal compressible MHD}\label{MHD}

Let us consider the $3D$ compressible ideal MHD equations. These are composed of
\bel{mhd1}
\rho\frac{\mathrm{D}\bu}{\mathrm{D}t} = (\mbox{curl}\,\bB) \times\bB - \bnabla P\qquad\mbox{and}\qquad 
\frac{\mathrm{D}\rho}{\mathrm{D}t} + \rho\mbox{div}\,\bu = 0\,,
\ee
together with the induction equation for the magnetic field $\bB$ (as well as the divergence-free condition $\mbox{div}\,\bB = 0$) and the equation for the specific internal energy $e$
\bel{mhd2}
\partial_t \bB = \mbox{curl}\, (\bu \times \bB) = \bB\cdot\bnabla\bu - \bB\divu - \bu \cdot \nabla \bB\,, \qquad\rho\frac{\mathrm{D}e}{\mathrm{D}t} = - P\mbox{div}\,\bu\,.
\ee 
The energy density $\mathcal{E}_{B}$ is an extension of (\ref{ed1}) to include $\bB$
\bel{mhd3}
\mathcal{E}_{B} = \shalf |\bB|^{2}+\mathcal{E}\,,
\ee
where $\mathcal{E} =  \rho\left(\shalf|\bu|^{2} + e\right)$ is the fluid energy density. Then we find that
\bel{mhd4}
\partial_{t}\mathcal{E}_{B} +
\mbox{div}\left\{\bu\left(\mathcal{E}_{B} + P\right) + \shalf |\bB|^{2} \bu - (\bu \cdot \bB) \bB \right\} = 0\,,
\ee
from which it follows that the full energy 
\bel{mhd5}
E_{0,B} = \I \left[\shalf |\bB|^{2}+\rho\left(\shalf|\bu|^{2}+e\right)\right]\,dV
\ee
is conserved. 
\par\smallskip%\noindent
The canonical cross-helicity $\mathcal{H}_{c} =\I \bu\cdot\bB\,dV$ is a pseudo-scalar \citep{Pouquet1976,Moffatt1978,Berger1984}. In parallel with (\ref{Hdef}), we introduce the following generalised cross-helicity 
\bel{magHdf}
\capHC = \I \rho\bu\cdot\bB\,dV\,,
\ee 
with the cross-helicity density defined as $h_{c} = \rho\bu\cdot\bB$. The MHD equivalent of the potential vorticity is $q_{c} = \bB\cdot\bnabla\rho$. One can check that the magnetic potential vorticity satisfies
\bel{qc}
\partial_t q_{c} + \Div \left( \bu q_c \right) + \Div \left( \rho \divu \bB \right) = 0\,.
\ee
The equation for $\rho\bu$ is
\bel{mhd4a}
\frac{\mathrm{D}(\rho\bu)}{\mathrm{D}t} + \rho\bu\mbox{div}\,\bu = (\mbox{curl}\bB)\times\bB - \bnabla P\,,
\ee
from which we deduce that
\beq{mhd4b}
\partial_{t}h_{c} &=& -\bB\cdot\left\{\bB\times \mbox{curl}\bB + \bnabla P + \rho\bu\mbox{div}\,\bu + \bu\cdot \bnabla(\rho\bu)\right\}\nonumber\\
&+& \rho\bu\cdot\left\{\bB\cdot\bnabla\bu - \bu\cdot\bnabla\bB - \bB\divu\right\}\\
&=& - 2h_{c}\divu - \bu\cdot\bnabla h_{c} - \bB\cdot\bnabla P +  \rho\bB\cdot\bnabla\left(\shalf|\bu|^{2}\right)\,.
\eeq
After re-arrangement we find 
\beq{mhd4d}
\partial_{t}h_{c} + h_{c}\divu +
\mbox{div}\bdJ_{c} &=& -\shalf q_{c}|\bu|^{2}\,.
\eeq
The equivalent of (\ref{Jhdef}) is 
\bel{lem3A}
\bdJ_{c} = h_{c}\bu + P\bB - \shalf\rho \bB |\bu|^{2}\,,
\ee
so with the definition $\sigma_{c} = -\shalf q_{c}|\bu|^{2} - h_{c}\divu$ we find
\bel{lem3B}
\partial_{t}h_{c} + {\rm div}\bdJ_{c} = \sigma_{c}\,.
\ee
This integrates to 
\bel{lem3C}
\frac{\mathrm{d}\capHC}{\mathrm{d}t} + \I h_{c}{\rm div}\bu\,dV = -\shalf\I q_{c}|\bu|^{2}dV\,.
\ee
As before, we note that the sign of the potential vorticity impacts the evolution of the cross helicity $\capHC$. Moreover, in regions of expansion $\capHC$ is decreasing, while it is increasing in regions of compression. Equation (\ref{lem3B}) is the equivalent of (\ref{lem2}). On the left-hand side, $\mbox{div}\,\bu$ weights the cross-helicity density. In the case of the inhomogeneous incompressible MHD equations, one finds a similar relation to (\ref{lem1}) and deduces that $q_{c}$ is a material constant. 

%%%%%%%%%%%%%%%%%
\section{\large Comments and conclusion}\label{comm}

\begin{table}
{\footnotesize
%  \begin{center}\def~{\hphantom{0}}
 \begin{adjustwidth}{0.0cm}{}
  \begin{tabular}{l|l|l|l|l}
             & $\mathcal{E}$   &   $h$ & $\bdJ$ & $\sigma$\\ \hline 
       Baro-Euler  & $\mathcal{E}_{0}$ & $h=\bu\cdot\bom$ & $\bdJ_{\pi} = h\bu + \bom\left(\Pi - \shalf|\bu|^{2}\right)$ & $0$\\ 
       II-Euler   & $\mathcal{E}_{0}$ & $h_{\rho}=\rho\bu\cdot\rho\bom$ & $\bdJ_{\rho} = h_{\rho}\bu + P\mbox{curl}\,(\rho\bu)-\shalf\bom\rho^{2}|\bu|^{2}$ & $\sigma_{\rho}=-q\rho|\bu|^{2}$\\
       Comp-Euler  & $\mathcal{E}$ & $h_{\rho}=\rho\bu\cdot\rho\bom$ & $\bdJ_{\rho} = h_{\rho}\bu + P\mbox{curl}\,(\rho\bu)-\shalf\bom\rho^{2}|\bu|^{2}$ & $\tilde{\sigma}_{\rho} = \sigma_{\rho} -2h_{\rho}\divu$\\
       MHD   & $\mathcal{E}+\shalf|\bB|^{2}$ & $h_{c}=\rho\bu\cdot\bB$ & $\bdJ_{c} = h_{c}\bu + \bB\left(P - \shalf\rho|\bu|^{2}\right)$ & $\sigma_{c} = -\shalf q_{c}|\bu|^{2} - h_{c}\divu$
  \end{tabular}
  \end{adjustwidth}
  \caption{\footnotesize The entries in the table represent the entropy-type relations $\partial_{t}h + \mbox{div}\,\bdJ = \sigma$ for the four different cases, which are the barotropic compressible Euler equations (Baro-Euler), the inhomogeneous incompressible Euler equations (II-Euler), the fully compressible Euler equations (Comp-Euler) and the ideal compressible MHD equations. Note that $\mathcal{E}_{0} = \shalf\rho|\bu|^{2}$ and $\mathcal{E} = \rho\left(\shalf|\bu|^{2}+e\right)$. }\label{tab:sum}
}
\end{table}
Before making some comments, let us summarise what we have found so far. Until now, most results on helicity dynamics for compressible flows have required the barotropic approximation. The main thread in this paper has been the investigation of a different definition of helicity density for compressible flows which takes the form $h_{\rho}=(\rho\bu)\cdot(\rho\bom)$ for which the barotropic condition is no longer required. The entries in Table \ref{tab:sum} show how the four cases are related. The inclusion of the $\rho^{2}$-term is crucial, as the results do not hold for other powers of $\rho$ in the helicity density. Moreover, we repeat that the evolution of $\capH$ (and $\capH_{c}$) can only be analysed provided there exist time intervals on which sufficiently smooth solutions exist for  either the compressible Euler equations (including the internal energy), the equations of compressible ideal MHD or the inhomogeneous incompressible Euler equations respectively. 
\par\smallskip%\noindent
The inclusion of the specific internal energy per unit mass $e$, while it does not appear explicitly in the calculations for the dynamics of $h_{\rho}$, impacts not only the $\mathcal{E}_{0}$-term through its appearance in the total energy equation (\ref{ed1}), but also the source term $\sigma_{\rho}$ defined in (\ref{Jhdef}). For the compressible Euler equations as given in equation \eqref{ce1a} to be fully determined, an equation of state for specifying the pressure $P$ is required. However, our results in this paper are independent of the choice of such an equation of state, which is why the dependence of the pressure on the density, temperature and entropy has been left unspecified. A typical non-barotropic choice might be the use of the ideal gas law, where $P$ would be a function of both density and temperature \citep{UM2005}.
\par\smallskip%\noindent
In the case of the inhomogeneous incompressible Euler equations, the new definition of $h_{\rho}$ and its behaviour leads us to introduce a resolution length scale $\lambda_{H}$ in (\ref{lamdefA}). This new length scale $\lambda_{H}$ is bounded from below in equation (\ref{lamest1B}) and it suggests a typical length scale on which helicity (and hence topological) variations occur. The two main features of the dynamics are the predominance of the potential vorticity $q$ which remains bounded in $L^{\infty}$, and in the fully compressible case the (average) sign of $\divu$. 

%%%%%%%%%%%%%%%%%%%%%%%%
\section*{\large Acknowledgements}
The authors are grateful to Keith Moffatt (Cambridge), Rahul Pandit (IISc Bangalore), Ian Roulstone (Surrey) and Dario Vincenzi (Nice) for helpful and encouraging discussions on this subject. The authors thank Robert M. Kerr (Warwick) for permission to use the figure of the helicity of a trefoil knot in the graphical abstract. They would also like to thank the Isaac Newton Institute for Mathematical Sciences, Cambridge, for support and hospitality during the programme ``Anti-diffusive dynamics\,: from sub-cellular to astrophysical scales'' where work on this paper was undertaken. This work was supported by EPSRC grant no EP/R014604/1. D.W.B. acknowledges support from the Cambridge Trust and the Cantab Capital Institute for Mathematics of Information.

%%%%%%%%%%%%%%%%%%%%%%%%
%%%%%%%%%%%%%%%%%%%
\bibliographystyle{jfm}
%\bibliography{jfm2esam}

\begin{thebibliography}{99}\itemsep 0mm
\small 

\expandafter\ifx\csname natexlab\endcsname\relax\def\natexlab#1{#1}\fi
\expandafter\ifx\csname selectlanguage\endcsname\relax
\def\selectlanguage#1{\relax}\fi

\bibitem[Berger and Field (1984)]{Berger1984} {\sc Berger, M. and  Field, G.} 1984 {The topological properties of magnetic helicity}, \textit{J. Fluid Mech.}, \textbf{147}, pp. 133--148.

\bibitem[Biferale {\em et al} (2012)]{Biferale2012} {\sc Biferale, L., Musacchio, S. and Toschi, F.} 2012 {Inverse energy cascade in three-dimensional isotropic turbulence}, \textit{Phys. Rev. Lett.}, \textbf{108}, pp. 164501.

\bibitem[Biferale and Titi (2013)]{Biferale2013} {\sc Biferale, L. and Titi, E. S.} 2013 {On the global regularity of a helical-decimated version of the 3D Navier-Stokes equations}, \textit{J. Stat. Phys.}, \textbf{151}, pp. 1089-1098.

\bibitem[Bjerknes (1898)]{VB1898} {\sc Bjerknes, V.} 1898 { \"{U}ber einen hydrodynamischen Fundamentalsatz und seine Anwendung besonders auf die Mechanik der Atmosph\"{a}re und des Weltmeeres}, \textit{Kongl. Sven. Vetensk. Akad. Handlingar}, \textbf{31}, pp. 1--35. 
% Atmosphäre

\bibitem[Boutros and Titi (2024)]{BT2024} {\sc Boutros, D. W. and Titi, E. S.} 2024 {On the conservation of helicity by weak solutions of the 3D Euler and inviscid MHD equations}, \textit{arXiv:2410.00813}.

\bibitem[Buaria {\em et al} (2022)]{BPB2022} {\sc Buaria, D., Pumir, P. and Bodenschatz, E.} 2022 {Generation of intense dissipation in high Reynolds number turbulence}, \textit{Phil. Trans. R. Soc. A} \textbf{380}, 20210088. 

\bibitem[Capocci {\em et al} (2023)]{Capocci2023} {\sc Capocci, D., Johnson, P. L., Oughton, S., Biferale, L. and Linkmann, M.} 2023 {New exact Betchov-like relation for the helicity flux in homogeneous turbulence}, \textit{J. Fluid Mech.}, \textbf{963}, pp. R1.

\bibitem[Cheskidov {\em et al} (2008)]{CCFS2008} {\sc Cheskidov, C., Constantin, P., Friedlander, S. and Shvydkoy, R.} 2008 {Energy conservation and Onsager's conjecture for the Euler equations}, \textit{Nonlinearity}, \textbf{21}, pp. 1233--1252.

\bibitem[Davies {\em et al} (2005)]{UM2005} {\sc Davies, T., Cullen, M. J. P., Malcolm, A. J., Mawson, M. H., Staniforth, A.,  White, A. A., and Wood, N.} 2005 {A new dynamical core for the Met Office's global and regional modelling of the atmosphere}, \textit{Quart. J. Royal Met. Soc.}, \textbf{131}, pp. 1759--1782.

\bibitem[Drivas and  Elgindi (2023)]{Drivas2023} {\sc Drivas, T. D. and  Elgindi, T. M.} 2023 {Singularity formation in the incompressible Euler equation in finite and infinite time}, \textit{EMS Surv. Math. Sci.}, \textbf{10}, pp. 1--100.

\bibitem[Elgindi (2021)]{Elgindi2021} {\sc Elgindi, T. M.} 2021 {Finite-time singularity formation for $C^{1,\alpha}$ solutions to the incompressible Euler equations on $\mathbb{R}^3$}, \textit{Ann. Math.}, \textbf{194}, pp. 647--727.

\bibitem[Enciso {\em et al} (2016)]{Enciso2016} {\sc Enciso, A., Peralta-Salas, D. and de Lizaur, F. T.} 2016 {Helicity is the only integral invariant of volume-preserving transformations}, \textit{Proc. Natl. Acad. Sci. USA}, \textbf{113}, pp. 2035-2040.

\bibitem[Ertel (1942)]{Ertel1942} {\sc Ertel, H.} 1942 {Ein neuer hydrodynamischer Wirbelsatz}, \textit{Meteorol. Z.}, \textbf{59}, pp. 271--281.

\bibitem[Gaffet (1985)]{Gaffet1985} {\sc Gaffet, B.} 1985 {On generalized vorticity-conservation laws}, \textit{J. Fluid Mech.}, \textbf{156}, pp. 141--149.

\bibitem[Gibbon and Holm (2012)]{Gibbon2012} {\sc Gibbon, J. D. and Holm, D. D.} 2012 {Quasiconservation laws for compressible 3D Navier-Stokes flow}, \textit{Phys. Rev. E}, \textbf{86}, 047301.

\bibitem[Helmholtz (1858)]{Helm1858} {\sc Helmholtz, H.} 1858 {\"{U}ber Integrale der Hydrodynamischen Gleichungen, welche den Wirbelbewegungen entsprechen}, \textit{J. Reine Angew. Math.}, \textbf{55}, pp. 25--55.

\bibitem[Iyer {\em et al} (2019)]{ISY2019} {\sc Iyer, K. P., Sreenivasan, K. R. and  Yeung, P. K.} 2019 {Circulation in High Reynolds Number Isotropic Turbulence is a Bifractal}, \textit{Phys. Rev. X}, \textbf{9}, 041006. 

\bibitem[Kerr (2018)]{Kerr2018} {\sc Kerr, R. M.} 2018 {Topology of interacting coiled vortex rings}, \textit{J. Fluid Mech.}, \textbf{854}, pp. R2.

\bibitem[Kerr (2023)]{Kerr2023} {\sc Kerr, R. M.} 2023 {Sensitivity of trefoil vortex knot reconnection to the initial vorticity profile}, \textit{Phys. Rev. Fluids}, \textbf{8}, 074701.

\bibitem[Linkmann (2018)]{Linkmann2018} {\sc Linkmann, M.} 2018 {Effects of helicity on dissipation in homogeneous box turbulence}, \textit{J. Fluid Mech.}, \textbf{856}, pp. 79–102.

\bibitem[Mobbs (1981)]{Mobbs1981} {\sc Mobbs, S.} 1981 {Some vorticity theorems and conservation laws for non-barotropic fluids}, \textit{J . Fluid Mech.}, \textbf{108}, pp. 475--483.

\bibitem[Moffatt (1969)]{Moffatt1969} {\sc Moffatt, H. K.} 1969 {The degree of knottedness of tangled vortex lines}, \textit{J. Fluid Mech.}, \textbf{35}, pp. 117--129.

\bibitem[Moffatt (1978)]{Moffatt1978} {\sc Moffatt, H. K.} 1978 \textit{Magnetic Field Generation in Electrically Conducting Fluids}, Cambridge University Press.

\bibitem[Moffatt (1985)]{Moffatt1985} {\sc Moffatt, H. K. } 1985 {Magnetostatic equilibria and analogous Euler flows of arbitrarily complex topology, Part 1. Fundamentals}, \textit{J. Fluid Mech.}, \textbf{159}, pp. 359--378.

\bibitem[Moffatt (1990)]{Moffatt1990} {\sc Moffatt, H. K.} 1990 {The energy spectrum of knots and links}, \textit{Nature}, \textbf{347}, pp. 367--369.

\bibitem[Moffatt and Ricca (1992)]{MoffattRicca1992} {\sc Moffatt, H. K. and  Ricca, R.} 1992 {Helicity and the C{\u{a}}lug{\u{a}}reanu Invariant}, \textit{Proc. R. Soc. Lond. A.}, \textbf{439}, pp. 411--429.

\bibitem[Moffatt (2014)]{Moffatt2014} {\sc Moffatt, H. K.} 2014 {Helicity and singular structures in fluid dynamics}, \textit{Proc. Natl. Acad. Sci. USA}, \textbf{111}, pp. 3663--3670.

\bibitem[Moreau (1961)]{Moreau1961} {\sc Moreau, J. J.} 1961 {Constantes d'un \^{i}lot tourbillonnaire en fluid parfait barotrope}, \textit{C. R . Acad, Sci. Paris}, \textbf{252}, pp. 2810--2812.

\bibitem[Poincar\'e (1893)]{Poincare1893} {\sc Poincar\'e, H.} 1893 {Th{\'e}orie des tourbillons\,: le{\c{c}}ons profess{\'e}es pendant le deuxi{\`e}me semestre 1891-1892 (Vol. 11)}, \textit{Gauthier-Villars, Article 158}.

\bibitem[Pouquet {\em et al} (1976)]{Pouquet1976} {\sc Pouquet, A., Frisch, U. and L\'eorat, J.} 1976 {Strong MHD helical turbulence and the nonlinear dynamo effect}, \textit{J. Fluid Mech.}, \textbf{77}, pp. 321--354.

\bibitem[Rogachevskii (2021)]{IR2021} {\sc Rogachevskii, I.} 2021 \textit{Introduction to Turbulent Transport of Particles, Temperature and Magnetic Fields}, Cambridge University Press. 

\bibitem[Roulstone and Norbury (2013)]{RN2013} {\sc Roulstone, I. and  Norbury, J.} 2013 \textit{Invisible in the Storm\,: The Role of Mathematics in Understanding Weather}, Princeton University Press.

\bibitem[Salmon (1988)]{Salmon1988} {\sc Salmon, R.} 1988 {Hamiltonian fluid mechanics}, \textit{Annu. Rev. Fluid Mech.}, \textbf{20}, pp. 225--256.

\bibitem[Sideris (1985)]{Sideris1985} {\sc Sideris, T.} 1985 {Formation of singularities in three-dimensional compressible fluids}, \textit{Comm. Math. Phys.}, \textbf{101}, pp. 475--485.

\bibitem[Silberstein (1896)]{Silber1896} {\sc Silberstein, L.} 1896 {O tworzeniu si{\k{e}} wir{\'o}w w p{\l}ynie doskona{\l}ym} (On the creation of eddies in an ideal ﬂuid), \textit{Proc. Cracow Acad. Sci.}, \textbf{31}, pp. 325--335. Also published as: Über die Entstehung von Wirbelbewegungen in einer reibungslosen Flüssigkeit, \textit{Bull. Int. l’Acad. Sci. Cracovie}, pp. 280–290.

\bibitem[Thomson (1868)]{Kelvin1868} {\sc Thomson, W.} 1868 {On vortex motion}, \textit{Trans. Roy. Soc. Edin.}, \textbf{25}, pp. 217--260.

\bibitem[Thorpe {\em et al} (2003)]{Thorpe2003} {\sc Thorpe, A. J., Volkert, H. and  Ziemianski, M. J.} 2003 {The Bjerknes' circulation theorem\,: a historical perspective}, \textit{Bull. Amer. Met. Soc.}, April, pp. 471--480.

\bibitem[Vallis (2019)]{Vallis2019} {\sc Vallis, G. K.} 2019 \textit{Essentials of Atmospheric and Oceanic Dynamics}, Cambridge University Press.

\bibitem[Webb (2018)]{Webb2018} {\sc Webb, G.} 2018 \textit{Magnetohydrodynamics and Fluid Dynamics: Action Principles and Conservation Laws}, Springer.

\bibitem[Woltjer (1958a)]{Woltjer1958} {\sc Woltjer, L.} 1958 {A theorem on force-free magnetic fields}. \textit{Proc. Natl. Acad. Sci. USA}, \textbf{44}, pp. 489--491.

\bibitem[Woltjer (1958b)]{Woltjer1958b} {\sc Woltjer, L.} 1958 {On hydromagnetic equilibrium}. \textit{Proc. Natl. Acad. Sci. USA}, \textbf{44}, pp. 833--841.

\bibitem[Yahalom (2017)]{Yahalom2017} {\sc Yahalom, A.} 2017 {A conserved local cross-helicity for non-barotropic MHD}, \textit{Geo. Astrophys. Fluid Dyn.} \textbf{111:2}, pp.  131--137.

\bibitem[Yeung and Ravikumar (2020)]{PK2020} {\sc Yeung, P. K. and  Ravikumar, K.} 2020 {Advancing understanding of turbulence through extreme-scale computation\,: Intermittency and simulations at large problem sizes}, \textit{Phys. Rev. Fluids}, \textbf{5}, 110517.

\bibitem[Yokoi (2013)]{Yokoi2013} {\sc Yokoi, Y.} 2013 {Cross-helicity and related dynamo}, \textit{Geo. Astrophys. Fluid Dyn.}, \textbf{107}, pp. 114--184.

\end{thebibliography}

\end{document}